\documentclass[11pt,fleqn,titlepage]{article}
\usepackage[citation-order, initials]{amsrefs}   
\usepackage{latexsym,amsmath,amssymb,amsfonts,epsfig}
\usepackage[all]{xy}

\providecommand{\CC}{{\mathbb{C}}}
\providecommand{\RR}{{\mathbb{R}}}

\providecommand{\ZZ}{{\mathbb{Z}}}

\providecommand{\QED}{{\hfill{} $\Box$}}

\providecommand{\FF}{{\mathcal F}}

\providecommand{\KK}{{\mathcal K}}

\newtheorem{definition}{Definition}

\newtheorem{proposition}[definition]{Proposition}

\newtheorem{theorem}{Theorem}

\begin{document}

\title{Noncommutative Topology and the World's Simplest Index Theorem}
\author{Erik van Erp \\
Dartmouth College \\
erik.van.erp@dartmouth.edu}

\date{March 2010}
\maketitle

\section{Introduction: Hypoelliptic Fredholm Index Theory}

In Euclidean space, waves propagate freely in all directions.
A standing wave with frequency $\omega$ is a solution of the eigenvalue problem of Laplace's equation,
\[ -\frac{\partial^2u}{\partial x^2} - \frac{\partial^2u}{\partial y^2} -\frac{\partial^2u}{\partial z^2} = \omega^2 u.\]
For any value of $\omega$ the solution space of this PDE is infinite dimensional,
consisting of superpositions of plane waves traveling in arbitrary directions. 
The degrees of freedom become severely restricted when we place the wave equation on a {\em closed manifold}.
Suddenly, solutions only exist for a discrete set of eigen-frequencies $\omega$, and even with the correct choice of $\omega$ the solution space is always finite dimensional.
The reason is intuitively obvious: unlike the situation in $\RR^3$, waves ``wrap around'' a closed manifold, and a standing wave can only exist if it wraps around in just the right way, always remaining exactly in phase with itself.
This lucky coincidence is the exceptional  case, and generically the local propagation law (the PDE) is in conflict with the obstruction presented by the global topology.

The index theorem of Atiyah and Singer is the ultimate expression of this connection between the local features of the PDE and the global topology of the manifold.
The {\em index} of a partial differential operator $P$ is not exactly the dimension of its solution space.
One must correct this dimension by subtracting the number of linearly independent conditions that are required of the ``source'' $v$ if the inhomogeneous equation $Pu=v$ is to have solutions.
Thus, the analytic index of $P$ is defined as
\[  {\rm Index}\, P = {\rm dim\; Kernel}\,P - {\rm codim\; Range}\,P.\]
Atiyah and Singer found a cohomology class $[\sigma(P)]$ associated to the principal symbol (i.e., the highest order part) of $P$, and expressed the analytic index of $P$ in terms of this cohomology class by means of a purely topological formula.
This {\em topological index} can often be computed explicitly with the tools of algebraic topology, and finding a positive index may be the only way to prove that solutions to the PDE {\em exist}---which, as we have seen, is not the generic situation.

A limitation of the theorem of Atiyah and Singer is that it only applies to {\em elliptic} operators.
Recall that while Laplace's equation is elliptic, the heat equation is parabolic, and the wave equation hyperbolic. 
What Laplace's equation and the heat equation (but not the wave equation) have in common is that all solutions---even weak solutions in the sense of distribution theory---are smooth functions.
Both equations are {\em hypoelliptic}.
Lars H\"ormander pioneered a deep investigation of hypoelliptic equations. 
It is not true in general that every hypoelliptic operator on a closed manifold is a {\em Fredholm operator},
i.e., an operator with a well-defined and finite analytic index.
But the methods typically used to prove hypoellipticity also imply Fredholmness.
Given these facts, it seems natural to search for an index theorem along the lines of the Atiyah-Singer formula for important classes of hypoelliptic Fredholm operators that appear in the literature.

A positive indication that this is possible was obtained by H\"ormander himself.
In 1971, almost ten years after Atiyah and Singer published their result, H\"ormander showed that the formula of Atiyah and Singer applies to hypoelliptic operators of type $(\rho, \delta)$
with $0\le 1-\rho\le \delta < \rho\le 1$ \cite{Ho71}.
In local coordinates, the full symbols of such operators are invertible outside a compact set, and these symbols can be glued together, by a partition of unity, to a well-defined class in $K^0(T^*M)$ to which the Atiyah-Singer formula applies.

Despite a growing literature devoted to the analysis of hypoelliptic operators, for over thirty years no progress was made in hypoelliptic index theory since the theorem of H\"ormander.
Only about a decade ago did Charles Epstein and Richard Melrose re-open the investigation \citelist{\cite{Me97} \cite{Ep04}}.
Their goal was to derive an index formula for the index of certain hypoelliptic operators associated to {\em contact structures}. 
The context for their work was not H\"ormander's theorem, but Boutet de Monvel's  index theorem for {\em Toeplitz operators} on strictly pseudoconvex CR manifolds \cite{Bo79}. 
Boutet de Monvel's theorem is essentially the Atiyah-Singer formula, except that the symbol of an elliptic operator is replaced by an expression derived from the symbolic calculus of the Toeplitz algebra.

Epstein and Melrose's goal was to extend Boutet de Monvel's formula to the {\em Heisenberg algebra} on a contact manifold. 
The Heisenberg algebra is a nonstandard pseudodifferential theory that extends the Toeplitz algebra and (like any pseudodifferential theory) includes also the algebra of differential operators.
The work of Epstein and Melrose (unfortunately largely unpublished) is a veritable tour de force, and they made considerable progress towards a solution.
Their final formula contains several factors that are not present in the formula of Atiyah and Singer, and that are very hard to compute in practice.
One of the key differences between the Atiyah-Singer and Boutet de Monvel formulas and the formula of Epstein-Melrose is that the latter is applied to a differential form (representing the operator) that is {\em not closed}, and therefore not a proper cohomology class.
We suspect that this is the primary reason for the appearance of the various mysterious correction factors.

In the present article we discuss a novel approach to the hypoelliptic index which has its roots in {\em noncommutative topology}.
As we will see, this modern perspective provides a satisfactory resolution of the problem researched by Epstein and Melrose, and shows that, after all, the Atiyah-Singer formula---just as it is, without bells and wistles---suffices to express the topological index,
just as in the index theorems of H\"ormander and Boutet de Monvel.
Moreover, our index theorem for contact manifolds appears as a special case of a single theorem in noncommutative topology.
Translating this noncommutative index theorem into classical algebraic topology can only be done on a case-by-case basis, resulting in theorems that no longer look alike
(a kind of ``symmetry breaking'' one could say).

We start our exposition with the statement of a new index formula for second order hypoelliptic PDEs on 3-manifolds.
This theorem contains no trace of its origin in noncommutative topology. 
To those who know elliptic theory,
the example will seem strange and unfamiliar,
and it will be hard to guess how this index theorem is related to the formula of Atiyah and Singer.

At the same time, the index formula presented here is as ``low-brow'' as we can make it, involving nothing more complicated than the notion of a winding number.
I am not aware of an index formula for {\em elliptic} differential operators that is this easy to state.\footnote{One could certainly make the case that the index formula for Toeplitz operators on a circle is even simpler. 
We maintain that differential operators are `simpler' than Toeplitz operators, if for no other reason than to justify the title of the article.
}
We hope this convinces the reader that, just because the proof involves noncommutative topology, the result is not necessarily esoteric.  

In the course of this article we gradually step back from this example, one step at a time.
The first step is to understand why the specific operators in our example are Fredholm. 
Here we discuss some classical results from analysis on nilpotent groups.
Our treatment of the {\em Heisenberg calculus} here is perhaps somewhat idiosyncratic, because we are mainly interested in applications to index problems as opposed to hypoelliptity.
This makes the exposition easier, and it suffices for our purposes.

Once we have a better understanding of the Fredholm theory of our operators, we can show that the winding numbers in our formula truly arise from an application of the Atiyah-Singer formula.
The insight that the Atiyah-Singer formula is all you need is precisely our main advance since the work of Epstein and Melrose.

The next section explains {\em why} the Atiyah-Singer formula works so universally,
and presents the noncommutative index theorem that is ultimately behind the simple formula involving winding numbers presented in the first section.
This ``mother of all index theorems'' says, in essence, that (1) there is only one topological index, namely the one calculated by Atiyah and Singer, and  (2) the nature of the cohomology class to which the formula is applied will depend on the details of the calculus.
Moreover, finding the right expression of this cohomology class in specific cases involves a precisely defined computational problem in noncommutative topology.
The level of difficulty of this computation depends on the spectral complexity of the algebra of symbols (we will explain this to some degree). 
The computation is easiest for elliptic operators (one recovers the class in $K^0(T^*M)$ defined by Atiyah and Singer), almost as easy for foliations (carried out in \cite{vEc}), and quite a bit harder for contact manifolds (Boutet de Monvel's theorem, generalized by Epstein-Melrose to the Hermite algebra, and now fully extended to the Heisenberg algebra).

Now that we have found the source of these theorems,
can we produce genuinely new ones?
We have explored this question in a couple of directions.
The structure theory of Type I $C^*$-algebras predicts the level of complexity of the index problem in individual cases.
We have explored some of the cases that are predicted to be most tractable.
In one direction, we found that Boutet de Monvel's Toeplitz index formula can be generalized to manifolds with {\em higher codimension contact structures},
such as appear at the boundary of quaternionic hyperbolic space. 
In another direction, our theory suggests that
{\em Heisenberg manifolds} should not be too much harder than contact manifolds.
The monograph [BG88] of Beals and Greiner is devoted to the analysis of second order hypoelliptic operators of a slightly more general class than the ones we study here.
This is work in progress. It seems that the index problem for Beals-Greiner operators may be just tractable, but our preliminary results are sufficiently complicated to discourage us from trying to push the method any further.

\section{The World's Simplest Index Theorem}

We start our discussion with a class of second order PDEs on 3-manifolds.
Let $M$ be a closed 3-manifold.
Consider a scalar second order differential operator $P$ on $M$,
represented in local coordinates as
\[ P = \sum_{i,j=1}^3 a_{ij} \frac{\partial^2}{\partial x_i \partial x_j} +
\sum_{i=1}^3 b_{i} \frac{\partial}{\partial x_i} + c.\]
The coefficients $a_{ij}, b_i, c$ are smooth functions on $M$.
If $(a_{ij})$ is a real symmetric matrix with  negative eigenvalues---so that $P$ is essentially a Laplacian---then the operator $P$ is elliptic and hence Fredholm.
However, the Fredholm index
\[ {\rm Index}\,P = {\rm dim\;Ker}\,P - {\rm dim\;Coker}\,P\]
of a Laplacian is always zero and so the index theory is not very interesting.

Let's see what happens if the matrix $(a_{ij})$ is degenerate, with two negative eigenvalues and one equal to zero.
We can cover $M$ by open sets in which $P$ is represented as
\[ P = -X^2-Y^2+{\rm lower\, order\, terms}.\]
Here $X$ and $Y$ are local vector fields on $M$ corresponding to eigenvectors of $(a_{ij})$.
Clearly, the operator $P$ is not elliptic. 
But a theorem of H\"ormander \cite{Ho67} provides necessary and sufficient conditions under which a second order differential operator (with real coefficients) is hypoelliptic. 
If the bracket $Z=[X,Y]$ is linearly independent of $X, Y$ then the role of the missing $-Z^2$ derivative is in some sense already taken care of by $-X^2-Y^2$.

\begin{proposition}\label{prop:bracket}
If the local vector fields $X, Y$ together with their bracket $[X,Y]$ span $TM$
(in each open set for which $P$ has a presentation as above),
and if $P$ has real coefficients,
then $P$ is hypoelliptic and Fredholm.
\end{proposition}
But operators that satisfy H\"ormander's condition also have zero index 
(as follows, for example, from the index formula we present below).
To get hypoelliptic Fredholm operators with a nonzero index, a more sophisticated analysis is needed, and we must include the {\em lower order terms} in our discussion.
So let us assume that $P$ is locally of the form
\[ P = -X^2-Y^2+i\alpha X + i\beta Y + i\gamma Z + \delta,\]
where we always have $Z=[X,Y]$, and require that $X, Y, Z$ span $TM$ in their domain of definition.
We allow {\em complex} values for the coefficients,
and therefore H\"ormander's bracket condition alone is not sufficient to guarantee Fredholmness.
Also remember that we do {\em not} require that the vector fields $X, Y, Z$ can be chosen {\em globally}.
Both of these facts will turn out to be crucial for constructing examples of operators with non-zero index.

As is easy to verify, the 2-plane bundle $\xi$ spanned by the local vector-fields $X, Y$ is a well-defined global vector bundle (it is dual to the characteristics of $P$).
Moreover, H\"ormander's bracket condition for local sections in $\xi$ is equivalent to $\xi$ being a {\em contact structure} on $M$.
By a foundational result in contact topology, every closed orientable 3-manifold admits a contact structure (Martinet's theorem).
In fact, there exist contact structures $\xi$ in every homotopy type of oriented 2-plane bundles in $TM$ (Lutz's theorem).
Thus, operators $P$ of the type we are interested in exist on every orientable 3-manifold.

Orientability of $M$ is a necessary condition, because the triple $X, Y, Z$ provides a well-defined orientation for $M$
(it has the same orientation as $Y, X, [Y,X]=-Z$).
In what follows we make the unnecesary but simplifying assumption that $\xi$ is (globally) oriented as well.

We now make an elementary but surprising observation.

\begin{proposition}
The $Z$-coefficient $\gamma$ in the local presentations of the operator $P$ is a well-defined global function $\gamma\colon M\to \CC$. 
\end{proposition} 
{\bf Proof.}
Suppose that in some open set in $M$ we have two alternative presentations of $P$ as
\begin{align*}
 P &= -X^2-Y^2+i\alpha X + i\beta Y + i\gamma Z + \delta,\\
&= -A^2-B^2+i\alpha' A + i\beta' B + i\gamma' C + \delta,
\end{align*}
where $Z=[X,Y]$ and $C=[A,B]$, and where both $X, Y$ and $A, B$ are positively oriented.
The second order terms in the two representations can only agree if
the pair $A, B$ is obtained from $X, Y$ by an orthogonal transformation
\[ A=aX+bY,\; B=cX+dY,\]
with 
\[ \left( \begin{array}{cc} a&b\\c&d \end{array} \right) \in SO(2).\]
Direct calculation shows that
\begin{align*}
A^2+B^2 &= (a^2+b^2)X^2 + (ab+cd)(XY+YX) + (c^2+d^2)Y^2\\
&+(aX.b+cX.d)Y+(bY.a+dY.c)X\\
&=X^2+Y^2+\xi-{\rm terms}.
\end{align*}
and similarly
\begin{align*}
C = AB-BA &= (ad-bc)(XY-YX)+(aX.c-cX.a+bY.c-dY.a)X \\
&+ (aX.d-cX.b+bY.d-dY.b)Y \\
&= Z + \xi-{\rm terms}.
\end{align*} 
Combining these results we see that $\gamma=\gamma'$.

\QED

To those accustomed to the Fredholm theory of {\em elliptic} operators, the following theorem will sound surprising.

\begin{proposition}\label{prop:odd}
The operator $P$ is Fredholm if the range of the $Z$-coefficient $\gamma$ does not contain any odd integers.
\end{proposition}
The calculus behind Proposition \ref{prop:odd} further implies---not surprisingly once we accept the proposition---that the index of the operator $P$ depends only on the contact structure $\xi$ and the homotopy type of the function
\[ \gamma \;\colon\; M\to \CC\setminus \{\ldots, -5, -3, -1, 1, 3, 5, \ldots\}.\] 
These curious facts have been known since the late 1970s (see, for example, \cite{FS74, Ro78}).
Beals and Greiner \cite{BG88} devoted an entire monograph to the analysis of exactly these operators (in a slightly more general setting).
But until the work of Epstein and Melrose the extensive literature on hypoelliptic operators contains no mention of the following question: {\it Can we express the Fredholm index of $P$ as a  ``topological index'', in terms of the homotopy type of the contact structure $\xi$ and the map $\gamma$?} 
To explain the gap of two or three decades between the analysis of these operators and the resolution of the associated index problem,
recall that all {\em scalar elliptic} differential operators have zero index,
as do {\em all elliptic} differential operators (acting on sections in a vector bundle) on odd dimensional manifolds.
What we have here is a scalar operator on an odd dimensional manifold.
If there was some obvious way to reduce the index problem for these hypoelliptic operators to an elliptic problem, one would expect to find a zero index in all cases.
But now that the problem is solved we know that every 3-manifold (except for a homology sphere) admits hypoelliptic scalar differential operators of the type we are studying here that have non-trivial index.

The resolution of the problem can be stated very simply as follows.
(A small detail needed to understand the theorem: the Euler class of an oriented contact 2-plane bundle is always {\em even}.) 

\begin{theorem}\label{thm:simple}

Let $L$ be an oriented link in $M$ such that the 1-cycle $2[L]\in H_1(M,\ZZ)$ represents the Poincar\'e dual of the Euler class $e(\xi)\in H^2(M,\ZZ)$.
For each odd integer $k$ the finite collection of loops $\gamma\colon L\to \CC\setminus \{k\}$ has a winding number
\[ {\rm Ind}_{L,\gamma}\,(k) = \frac{1}{2\pi i}\int_{L}\,\frac{d\gamma}{\gamma-k}\;\in \ZZ.\]
The Fredholm index of $P$ is a $\ZZ$-linear combination of these winding numbers, \[ {\rm Index}\, P = \sum_{k \;\rm {odd}}\; k\,\cdot\,{\rm Ind}_{L,\gamma}\,(k).\]
\end{theorem}
Observe that if all coefficients of $P$ are real valued (as in Proposition \ref{prop:bracket}),
then $\gamma$ is purely imaginary and the map $\gamma\colon M\to \CC\setminus \{{\rm odds}\}$ is contractible.
Theorem \ref{thm:simple} implies that ${\rm Index}\,P=0$.
This is not surprising, and can probably be proven by a direct homotopy from $P$ to a self-adjoint operator. 
However, Theorem \ref{thm:simple} also implies that if the vector fields $X$ and $Y$ in the presentation   
\[  P = -X^2-Y^2+i\alpha X + i\beta Y + i\gamma Z + \delta,\]
are {\em globally} defined,
then also ${\rm Index}\,P=0$, regardless of the homotopy type of $\gamma$.
I am not aware of an elementary proof of this corollary.

\section{Changing Orders to Suit Our Needs}

Before outlining the proof of Theorem \ref{thm:simple}
we must digress and explain the calculus behind Proposition \ref{prop:odd}.
Why is it that Fredholmness of $P$ depends on the {\em lower order term} $i\gamma Z$? 
The answer, in a word, is that $i\gamma Z$ is included in the {\em highest order part} of $P$, if only we re-define ``order'' in a suitable way.  
(References for the calculus described in this section are \cite{Ta84, BG88, Ep04}). 
However, the perspective taken here is slightly different, and influenced by our work on index theory.)

The highest order part of a differential operator like
\[ P = \sum_{i,j=1}^3 a_{ij} \frac{\partial^2}{\partial x_i \partial x_j} +
\sum_{i=1}^3 b_{i} \frac{\partial}{\partial x_i} + c\]
at a point $m\in M$,
\[ P_m = \sum_{i,j=1}^3 a_{ij}(m) \frac{\partial^2}{\partial x_i \partial x_j},\]
is not well-defined as a differential operator on $M$.
The algebra of differential operators is only filtered, not graded.
Nevertheless, the highest order part of $P$ is well-defined and is usually taken to be a constant coefficient operator on the tangent fiber $T_mM\cong\RR^3$, independent of coordinate choices. 

All we need to do to include $i\gamma Z$ in the highest order part of our hypoelliptic operators is to change the filtration on the algebra of differential operators.
However, in so doing we also must reinvestigate what {\em kind} of object the highest order part of an operator is.

The {\em Heisenberg calculus} is, in essence, nothing more than the working out of the consequences of an alternative filtration on the algebra of differential operators on $M$.
The filtration is defined as follows: all vector fields in the direction of the contact field $\xi$ (like $X$ and $Y$) have $\xi$-order one, as usual, but {\em any} vector field that is not everywhere tangent to $\xi$ (like $Z$) has $\xi$-order {\em two}. 
The highest order part of $P$ in the Heisenberg calculus will then clearly include the $i\gamma Z$-term.

But we must also investigate what ``highest order part'' really means in this new calculus.
Abstractly, for any filtered algebra, the notion of ``highest order part'' refers to an element in the associated graded algebra.
So we need to find an analytic model for the associated graded algebra of the algebra of differential operators with the Heisenberg filtration.

Observe that in the associated graded algebra, smooth functions $f$ commute with all vector fields $A$, because the commutator $[A,f]=A.f$ is of order zero. 
This implies that elements in the graded algebra can be {\em localized} at points $m\in M$.
The result of this localization is, formally,
\[ P_m = -X(m)^2-Y(m)^2 + i\gamma(m) Z(m).\]
As before, the highest order part of $P$ is not an operator on $M$,
but it has a natural interpretation as a smooth family $P_m, m\in M$ of operators in the tangent fiber $T_mM$. 
The difference is that in the Heisenberg calculus the operator $P_m$ is not  a constant coefficient operator, but an operator that is {\em left invariant} for a nilpotent group structure on the tangent fiber.
One identifies $T_mM$ (or, more precisely, $\xi_m\oplus T_mM/\xi_m$) with the {\em Heisenberg group}, by imposing the commutation relations
\[ [X, Y] = Z,\; [X, Z] = [Y, Z] = 0.\]
Then the $\xi$-highest order part of $P$ is the smooth family $\{P_m, m\in M\}$ of left invariant homogeneous operators on the (graded) Heisenberg group 
\[ G_m=\xi_m\oplus T_mM/\xi_m \cong T_mM.\]
If one wants a more explicit formula for the model operator $P_m$, one could substitute
\[
 X(m) = \frac{\partial}{\partial x} - \frac{1}{2} y\frac{\partial}{\partial z},\;
 Y(m) = \frac{\partial}{\partial y} + \frac{1}{2} x\frac{\partial}{\partial z},\;
 Z(m) = \frac{\partial}{\partial z}.
\]
While this substitution is common in the analytic literature, 
it is not useful for our purposes.
We prefer to retain the simpler algebraic expression of $P_m$ in terms of $X, Y, Z$.

We can now state the main results about hypoelliptic operators in the Heisenberg calculus.
\begin{definition}
A $\xi$-homogeneous invariant operator on a graded nilpotent group is a {\em Rockland operator} if $\pi(P_m)$ is invertible for all irreducible unitary representations $\pi$ of the group, except for the trivial representation.
\end{definition}
Observe that for the abelian group $\RR^n$ unitary representation theory amounts to Fourier theory, and a Rockland operator on $\RR^n$ is just a homogeneous elliptic constant coefficient operator.  
In this language, one could define ellipticity by saying that $P$ is elliptic if each of the model operators $P_m$ in the highest order part of $P$ is a Rockland operator on the abelian group $T_mM$.
This generalizes beautifully to the Heisenberg calculus.
\begin{definition}\label{def:hypo}
A differential operator $P$  is {\em $\xi$-elliptic} if all the model operators $P_m$ in its $\xi$-highest-order in the Heisenberg calculus are Rockland operators.
\end{definition}

\begin{proposition}\label{prop:Rockland}
Let $(M, \xi)$ be a closed contact manifold, and $P$ a scalar differential operator on $M$.
If for all $m\in M$ the $\xi$-highest order part $P_m$ is a Rockland operator, 
then $P$ is a hypoelliptic Fredholm operator.
\end{proposition} 
Let's apply this to our example:
When is $P_m=-X^2-Y^2+i\gamma(m)Z$ a Rockland operator?
We need explicit formulas for the representations of the Heisenberg group.
Since we are interested in the representation of the differential operator $P_m$,
it is most convenient to express our formulas in terms of the Lie algebra. 
The Heisenberg group has two families of irreducible representations. 
First, there is an infinite family of scalar representations
with two continuous parameters $(x, y)$, given by
\[ \pi_{(x,y)}(X) = ix,\;\pi_{(x,y)}(Y) = iy,\; \pi_{(x,y)}(Z) = 0.\]
When we apply these to $P_m$ we find 
\[ \pi_{(x,y)}(P_m) = x^2+y^2,\]
which is invertible unless $(x,y)=(0,0)$ (the trivial representation). 
We see that invertibility of the nontrivial scalar representations of $P_m$
amounts to ``partial ellipticity'' of $P$ in the $\xi$-directions,
and it is built into our definition of $\xi$ from $P$.

The second family of irreducible representations of the Heisenberg group are labeled by a single parameter $t\in \RR\setminus\{0\}$ ($t=0$ corresponds to the entire collection of scalar representations).
The representation space is the {\em Bargman-Fock} space.
It is the completion of the space of holomorphic polynomials in one complex variable $z\in \CC$. 
The monomials $z^q/q!$ form an orthonormal basis. 
The normalization is chosen such that the operators $z$ and $\partial/\partial z$ are adjoints. 

If $t>0$  the representation $\pi_t$ is defined by
\[ \pi_t(X) = i \,\sqrt{\frac{t}{2}}\, \left(z+\frac{\partial}{\partial z} \right),\;
\pi_t(Y) = \sqrt{\frac{t}{2}}\, \left(z-\frac{\partial}{\partial z} \right),\;
\pi_t(Z) = it,\]
while $\pi_{-t}$ is the conjugate representation of $\pi_t$.
We find for $t>0$,
\[ \pi_{t}(P_m) = -\pi_{t}(X)^2-\pi_{t}(Y)^2+i\gamma(m)\pi_{t}(Z) = t\left(2z\frac{\partial}{\partial z}+1-\gamma(m)\right).\]
As one easily verifies, the basis vectors $z^q/q!$ are eigenvectors for the operator $\pi_{t}(P_m)$ with eigenvalues $2q+1-\gamma(m)$.
It follows that $\pi_{t}(P_m)$ is invertible for all $t>0$ iff $\gamma(m)$ is not a positive odd integer.

Observe that, due to the {\em homogeneity} of $P$, we have
\[ \pi_t(P_m) = t \; \pi_{+1}(P_m), t>0,\]
as can be seen explicitly in the formula above. Likewise,
\[ \pi_{-t}(P_m) = t\; \pi_{-1}(P_m), t>0.\]
For this reason, it suffices to verify invertibility of $\pi_{+1}(P_m)$ and $\pi_{-1}(P_m)$. 
A convenient way of dealing with the conjugate representation $\pi_{-1}$ is as follows.
Consider the anti-automorphism of the Lie algebra of the Heisenberg group given by
\[ X^{op} = X,\; Y^{op}=Y, Z^{op}=-Z.\]
This induces an anti-automorphism of the algebra of invariant operators on $G$. 
Therefore invertibility of $\pi_{-1}(P_m)$ is equivalent to invertibility of $\pi_{+1}(P_m^{op})$. Thus, all we need to do is verify that $\pi_{+1}(P_m)$ and $\pi_{+1}(P_m^{op})$ are invertible.
The advantage of this description is that we can work with the single representation $\pi_{+1}$.
This will also prove useful when we get to our index formula. 

To finish the discussion, simply reversing the sign of $\gamma$ in $\pi_{+1}(P_m)$ we obtain,
\[ \pi_{+1}(P_m^{op})=2z\frac{\partial}{\partial z}+1+\gamma(m).\]
This operator is invertible iff $\gamma(m)$ is not a {\em negative} odd integer.  
So, by Definition \ref{def:hypo}, $P$ is $\xi$-elliptic precisely if $\gamma(m)$ is not an odd integer for any $m\in M$.
Then Proposition \ref{prop:odd}, which seemed bizarre from the perspective of elliptic theory, appears as a special case of Proposition \ref{prop:Rockland}.
Once appropriately generalized, what seemed strangely different is seen to be a perfect analogy.

\section{The Atiyah-Singer Formula (Sort Of)}

The previous discussion suggests that our hypoelliptic index problem should be rephrased, in the context of the Heisenberg calculus, as follows:
{\em Find a topological index that expresses the Fredholm index of $P$ in terms of its $\xi$-highest order part}
\[ \sigma_\xi(P)=\{P_m, m\in M\}.\]
This was, in essence, what Epstein and Melrose set out to do.
Contrary to what is suggested by the complicated formula they found (announced in \cite{Me97}), it turns out that the topological index for $\xi$-elliptic operators is no different from the topological index computed by Atiyah and Singer.
\begin{theorem}[\cite{vEb}]\label{thm:contact}
Let $P$ be a $\xi$-elliptic differential operator on a closed contact manifold $(M, \xi)$. Then
\[{\rm Index}\, P = \int_M\,{\rm Ch}(\sigma_\xi(P)) \wedge {\rm Td}\,(M),\]
where
\[ [\sigma_\xi(P)] =   \left[\frac{\pi(P_m)}{\pi(P_m^{op})}\right] \,\in K^1(M).\]
\end{theorem}
This is the well-known cohomological formula of Atiyah and Singer,
except that the principal symbol $\sigma(P)$ of an elliptic operator is replaced by the Heisenberg symbol $\sigma_\xi(P)$.

To explain our formula we need to clarify the meaning of the quotient
\[ [\pi(P_m)\pi(P_m^{op})^{-1}] \in K^1(M)\] 
which appears under the Chern character. 
Rather than explain this in general, we will illustrate its meaning by deriving Theorem \ref{thm:simple} from Theorem \ref{thm:contact} by explicitly computing the appropriate $K$-theory class.

So assume again that $P$ is a $\xi$-elliptic scalar operator on a closed contact $3$-manifold, locally presented as    
\[ P = -X^2-Y^2+i\alpha X + i\beta Y + i\gamma Z + \delta.\]
We have calculated $\pi_{+1}(P_m)$ and $\pi_{+1}(P_m^{op})$,  
and we see that the quotient of these two operators is a diagonal operator with eigenvalues
\[ \frac{2q+1-\gamma(m)}{2q+1+\gamma(m)}.\]
These scalars act on the basis vector $z^q/q!$, which is canonically interpreted as a section in the $q$-th tensor power $(\xi^{1,0})^{\otimes q}$.
Here $\xi^{1,0}$ denotes the oriented 2-plane bundle $\xi$ thought of as a complex line bundle.

Recall that, in general, a cocycle in $K^1(M)$ is represented by a pair $[E, \sigma]$ consisting of a vectorbundle $E$ on $M$ and an automorphism $\sigma$ of $E$.  
In the case of our $\xi$-elliptic operator $P$,
we find, 
\[ \left[\frac{\pi(P_m)}{\pi(P_m^{op})}\right] = \sum_{q=0}^\infty\, 
               \left[\left(\xi^{1,0}\right)^{\otimes q}, \frac{2q+1-\gamma}{2q+1+\gamma}\right] \;\in K^1(M). \]
This is the meaning of this expression as it appears in Theorem \ref{thm:contact}.
Observe that the right hand side only {\em appears} to be an infinite sum,
because for sufficiently large $q$ the automorphisms are close to (and therefore homotopic to) the trivial automorphism. 

It is this computation of a $K$-theory class associated to $P$ that is unique to the Heisenberg calculus for contact manifolds.
Once we have the right $K$-theory class we can apply the formula of Atiyah and Singer. 
Using the fact that the Chern character is a ring homomorphism (from $K$-theory to deRham cohomology) we find
\[ {\rm Ch}(\sigma_\xi(P)) = \sum_{q=0}^\infty\, 
               {\rm Ch}\,(\xi^{1,0})^q\wedge {\rm Ch}\,\left(\frac{2q+1-\gamma}{2q+1+\gamma}\right).\]
Because $\xi^{1,0}$ is a complex line bundle on a 3-manifold we have
\[ {\rm Ch}\,(\xi^{1,0})^q = (1 + c_1(\xi^{1,0}))^q = 1 + q\cdot e(\xi),\]
where $e(\xi)$ denotes the Euler class of the oriented real 2-plane field $\xi$. The (odd) Chern character of an invertible function $u\colon M\to \CC$ is given by 
\[ {\rm Ch}\,(u) = \frac{1}{2 \pi i} d\log{u}.\]
Because ${\rm dim}\,M = 3$ we have ${\rm Td}\,(M) = {\rm Td}\,(\xi^{1,0})=1+e(\xi)/2$,
and the Atiyah-Singer formula gives
\begin{align*}
{\rm Index}\,P &= \int_M {\rm Ch}\,(\sigma_\xi(P)) \wedge {\rm Td}(M) \\
&= \sum_{q=0}^\infty\,\int_M  \frac{2q+1}{2}e(\xi) \wedge d\log{(2q+1-\gamma)} \,-\, \frac{2q+1}{2}e(\xi) \wedge d\log{(2q+1+\gamma)},
\end{align*}
which is equivalent to Theorem \ref{thm:simple}.\footnote{In a personal communication, Charles Epstein pointed out that the index formula for the special case of second order differential operators can also be derived by a deformation of the Heisenberg symbol of the differential operator to the symbol of a Toeplitz operator, to which Boutet de Monvel's index theorem applies.
It is not clear if or how this trick extends to the general case.
}

\section{A Factory for Hypoelliptic Index Theorems}

We have explained in what sense Theorem \ref{thm:simple} is an application of the topological index of Atiyah and Singer. 
From the point of view of classical analysis or topology it is very surprising that the Atiyah-Singer formula can be modified in such a straightforward manner.
Our next objective is to explain why this is so. 
As before, what appears strange from one point of view is clarified when seen from a higher level of abstraction.
In the present case, we must understand the problem in the context of {\em noncommutative topology}---specifically, the theory of $C^*$-algebras and analytic $K$-theory.
 
The idea of analytic $K$-theory is rooted in Atiyah's representation of $K^0(X)$ for compact $X$ as the set of homotopy classes of maps from $X$ into $\FF$,
\[ K^0(X) = [X, \FF],\]
where $\FF$ is the space of all Fredholm operators on a given (infinite dimensional, separable) Hilbert space.
Concretely, if $\{P_x, x\in X\}$ is a continuous family of Fredholm operators $P_x$ over $X$, then the family of kernels $E_x$ of $P_x$ defines a vector bundle $E$ over $X$ (possibly of nonconstant fiber dimension), and so does the family $F_x$ of cokernels. The ``index'' of this family can then be thought of as the formal difference 
\[ [E]-[F] \in K^0(X).\]
This ``family index'' is Atiyah's isomorphism.

The idea that $K$-theory elements can be represented by families of Fredholm  operators has been generalized and developed into analytic $K$-theory of $C^*$-algebras.
The limited space of this article does not allow us to explain these ideas here in any detail.
However, the example of Atiyah's theorem may help the reader accept the following facts without much explanation.
We have seen that the highest order part of an elliptic operator $P$ on $M$ is a family $\{P_m, m\in M\}$ of elliptic operators.
These operators are not Fredholm, because they act on functions $C_c^\infty(T_mM)$ on the Euclidean space $T_mM$, which is not compact.
But they define a ``Fredholm operator'' in a generalized sense.
First, one can think of the family as a {\em single} differential operator on $C_c^\infty(TM)$ which only differentiates in the direction of the fibers.
The {\em translation invariance} of the operators $P_m$ can be expressed algebraically by treating $C_c^\infty(TM)$ as a {\em convolution algebra}, where functions on $TM$ are multiplied by taking their convolution product in each fiber.
Then the family $\{P_m, m\in M\}$---taken as a single operator on $C_c^\infty(TM)$---commutes with the right module action of $C_c^\infty(TM)$ on itself.
Ellipticity of the individual operators $P_m$ then implies that this $C_c^\infty(TM)$-linear operator satisfies the axioms of a suitably generalized notion of ``Fredholm operator''.
One completes $C_c^\infty(TM)$ in a suitable norm to obtain the $C^*$-algebra $C^*(TM)$. 
The $C^*(TM)$-linear ``Fredholm operator'' $\{P_m, m\in M\}$ then defines an (unbounded) element in the analytic $K$-theory of $C^*(TM)$,
\[ [\{P_m, m\in M\}] \in K_0(C^*(TM)).\]
Of course, the Fourier transform in the fibers $T_mM$ of the tangent space gives an algebra isomorphism
\[ C^*(TM) \cong C_0(T^*M).\]
After Fourier transform the family of operators $\{P_m, m\in M\}$ is replaced by the function $\sigma(P)$ on $T^*M$ that is usually referred to as the ``principal symbol'',
and the analytic $K$-theory element
$[\{P_m, m\in M\}]$ corresponds exactly to the topological $K$-theory class $[\sigma(P)]$ defined by Atiyah and Singer,
\[ K_0(C^*(TM))\cong K^0(T^*M)\;\colon\; [\{P_m, m\in M\}] \mapsto [\sigma(P)].\]
What is the point of this translation?
After all, the definition of the topological class $[\sigma(P)]$ seems much simpler than that of the equivalent analytic class $[\{P_m, m\in M\}]$.
It surely is.
The point of the analytic construction is that it applies much more generally.
In particular, it applies without essential modification to the highest order part $\{P_m, m\in M\}$ of a $\xi$-elliptic operator $P$ on a contact manifold.
This time, the operators $P_m$ are translation invariant for a nilpotent group structure in the fibers of $TM$.
Let us denote by $T_\xi M$ the tangent space understood as the smooth family of Heisenberg groups $T_mM\cong \xi_m\oplus \RR$ that underly the Heisenberg calculus.
As before, we can form the convolution algebra $C_c^\infty(T_\xi M)$ and its completion $C^*(T_\xi M)$.
These algebras are noncommutative, but we don't have to worry about that.
The highest order part $\{P_m, m\in M\}$ of any operator $P$ in the Heisenberg calculus can be conceived as an operator on $C_c^\infty(T_\xi M)$ that is right-$C_c^\infty(T_\xi M)$ linear.
If, in addition, every $P_m$ is a Rockland operator (which makes $P$ $\xi$-elliptic), then one can prove that this family is indeed a ``Fredholm operator'' in the generalized sense of analytic $K$-theory. 
Therefore, as before,
\[ [\{P_m, m\in M\}] \in K_0(C^*(T_\xi M)).\]
So, at the very least, we have a $K$-theory element, albeit an analytically defined $K$-theory element as opposed to a topological one.
The goal is, then, to translate this into topology. 
But analytic $K$-theory is as much noncommutative {\em topology} as it is analysis.
In the present case, a noncommutative version of the Thom isomorphism in $K$-theory (the Connes-Thom isomorphism) implies that
\[ K_0(C^*(T_\xi M))\cong K^0(T^*M).\]
Therefore, at least theoretically, our analytically constructed $K$-theory class can be interpreted as a topological $K$-theory class,
\[ [\sigma_\xi(P)] \in K^0(T^*M).\]
In fact, the construction generalizes further.
Suppose that $M$ is equipped with an {\em arbitrary} sub-bundle $H\subseteq TM$, of arbitrary codimension. The bundle $H$ could be a foliation, or a contact structure, it could be related to a sub-Riemannian structure on $M$, or it could be a random bundle with no geometric significance whatsoever.
We can study hypoelliptic operators that are elliptic in the directions encoded by $H$, and of lower order in the transversal directions.
Such operators will be the elliptic elements in a generalized Heisenberg calculus associated to $H$ (where transversal vector fields are given order two.)
The highest order part for $H$-elliptic operators in this calculus is a family of invariant Rockland operators $P_m$ on nilpotent groups $G_m=H_m\oplus T_mM/H_m$,
and this family represents an analytic $K$-theory class 
\[ [\{P_m, m\in M\}] \in K_0(C^*(T_HM)).\]
Once again this group is naturally isomorphic to $K^0(T^*M)$, and so we have
\[ [\sigma_H(P)]\in K^0(T^*M).\]
Thus, using the tools of noncommutative topology, one can {\em always} construct a symbol class in $K^0(T^*M)$ for this type of hypoelliptic operators.
Then a very general index problem presents itself: {\em find the topological index 
\[ {\rm Ind}_H\;\colon\;K^0(T^*M)\to \ZZ\]
for $H$-elliptic operators.}
A priori, this topological index could depend on the geometric structure $H$.
But it doesn't. 
A beautiful argument involving deformation theory of $C^*$-algebras (this is where Lie groupoids play a key role) implies that, once we have identified the symbol of a hypoelliptic operator as a class in $K^0(T^*M)$, 
there is only one topological index.

\begin{theorem}\label{thm:topind}
Let $M$ be a closed manifold equipped with a subbundle $H\subseteq TM$.
Let $P$ be an $H$-elliptic operator on $M$---elliptic in the directions $H$, and of lower order transversally.
Then the highest order part $\sigma_H(P) = \{P_m, m\in M\}$ in the generalized Heisenberg calculus for $(M,H)$ defines a $K$-theory class
\[ [\sigma_H(P)] \in K_0(C^*(T_HM)) \cong K^0(T^*M)\] 
and the Fredholm index of $P$ is computed by the topological index of Atiyah and Singer,
\[ {\rm Index}\, P = \int_{T^*M} {\rm Ch}(\sigma_H(P)) \wedge {\rm Td}\,(M).\]
\end{theorem}
This theorem generalizes Boutet de Monvel's theorem for Toeplitz operators, and it is also the source of our proof of Theorem \ref{thm:contact}.
 Nevertheless, Theorem \ref{thm:topind}, by itself, is too abstract to be useful.
To turn it into a computable formula---say, something as concrete as Theorem \ref{thm:simple}---two obstacles need to be overcome.

The first obstacle is a purely analytic problem in the Heisenberg calculus.
Given an explicit differential operator $P$, 
it is not too hard to verify whether its $H$-symbol consists of Rockland operators $P_m$.
But that is not sufficient in a general setting.
One must still prove that if each $P_m$ is Rockland, then $P$ has an inverse in the generalized Heisenberg calculus.
This is known to be true if the manifold (with its structure $H$) can be locally identified with the graded nilpotent model group---such as is the case for foliations (Frobenius's theorem) and contact structures (Darboux's theorem).
But in general it can be quite tricky. 
How hard the problem is in general can be gleaned by skimming through the monograph of Beals and Greiner \cite{BG88}. 
The bulk of that publication is devoted to proving just this fact for second order $H$-elliptic operators in the case that $H$ is a subbundle of co-dimension one in $TM$.
To my knowledge, it is not known whether this fact (if all $P_m$ are Rockland then $P$ has an inverse in the calculus) is true in absolute generality.

The second obstacle is a problem in noncommutative topology.
In order to compute the Chern character of $\sigma_H(P)$ we must have an explicit understanding of the Connes-Thom isomorphism
\[ K_0(C^*(T_HM)) \cong K^0(T^*M).\]
Unfortunately, there is no general method for finding a topological representation of the analytically defined class.
Interestingly, the map from analytic to topological $K$-theory is, in essence, a generalized index---as in the example of Atiyah's analytic representation of elements in $K^0(X)$.
Therefore, this computation of the symbol class in $K^0(T^*M)$ is itself an index problem---not the index of $P$, but the index of its highest order part (which, as we have indicated, is a generalized Fredholm operator).
The tractability of this problem depends on the specifics of the geometric structure $(M,H)$, as we will see in our next and final section.

\section{Offspring: Foliations, Contact Manifolds, and More}

To derive Theorem \ref{thm:contact} from Theorem \ref{thm:topind} and thereby give a proof of the World's Simplest Index Theorem (Theorem \ref{thm:simple}),
we need to find an explicit formula for the isomorphism
\[ K_0(C^*(T_\xi M)) \cong K^0(T^*M)\]
if $\xi$ is a contact structure.
The computation depends on a careful analysis of the structure of the $C^*$-algebra $C^*(T_\xi M)$.

In the general case, if $H\subseteq TM$ is an arbitrary bundle, the structure of the $C^*$-algebra $C^*(T_HM)$ is clarified by the theory of {\em Type I} $C^*$-algebras.
A $C^*$-algebra is of Type I if it can be decomposed, by a series of successive short exact sequences, into commutative $C^*$-algebras, or algebras that are $K$-theoretically equivalent to commutative algebras.
This reduction to commutative algebras is precisely what establishes the connection between noncommutative and classical topology,
and will allow us to use topological $K$-theory to compute an analytic $K$-cocycle.
Thus, the structure theory of $C^*$-algebras tells us something significant about the geometric structure $(M,H)$, and about the difficulty of the related index problem.
It is reasonable to expect that the more terms there are in the decomposition series of $C^*(T_HM)$ (it is always finite), the more intractable the index problem will be.

In the trivial case where $H=TM$---the case of elliptic operators---the computation is fairly straightforward,
because of the isomorphism $C^*(TM)=C_0(T^*M)$ established by Fourier theory.
In general, $C^*(T_HM)$ is commutative if the groups $G_m$ are abelian, which is the case precisely if $[H,H]\subseteq H$, i.e., if $(M,H)$ is a {\em foliated manifold}. 
For foliations, the computation of the class $\sigma_H(P)\in K^0(T^*M)$ for $H$-elliptic operators (these are {\em not} elliptic!) is not much harder than it is for elliptic operators.
We have carried out this computation, and by this general methodology derive a hypoelliptic index theorem that is similar in spirit to the one proved by  H\"ormander for hypoelliptic operators of type $(\rho, \delta)$, except that hypoelliptic operators (and their parametrices) in the Heisenberg calculus are of type $(\frac{1}{2}, \frac{1}{2})$, a case not covered by H\"ormander's theorem (see \cite{vEc}).

From the prespective of the structure theory of $C^*(T_HM)$, contact structures are the next best thing, after foliations. 
If $H=\xi$ is a contact structure, there is a single short exact sequence that decomposes $C^*(T_\xi M)$, as follows,
\[ 0\to C_0(M\times \RR^{\times})\otimes \KK \to C^*(T_\xi M) \to C_0(\xi^*)\to 0.\]  
As before, we assume here that $\xi$ is co-oriented, i.e., $TM/\xi$ is a trivial line bundle.
Then $\RR^{\times}=\RR\setminus \{0\}$ refers to the punctured fiber of the dual of the line bundle $TM/\xi\cong M\times \RR$. 
The elementary algebra $\KK$ is the (abstract) $C^*$-algebra of compact operators on a separable Hilbert space.

To understand this sequence, it is helpful to focus on a single fiber of the bundle of algebras $C^*(T_\xi M)$. 
For every $m\in M$, the fiber is the group $C^*$-algebra $G^*(G_m)$ 
(the completion of the convolution algebra $C_c^\infty(G_m)$ in a suitable norm)
of the Heisenberg group $G_m=H_m\times \RR$.
This group algebra decomposes as
\[ 0\to C_0(\RR^*\setminus \{0\})\otimes \KK \to C^*(G_m) \stackrel{\phi}{\to} C_0(H_m^*) \to 0.\]
The quotient map $\phi$ is obtained by assembling all the scalar representations of the group $G_m$, parametrized by points in $H_m^*$.
The remaining irreducible representations of $G$ are labelled by the dual of the center $\RR\subset G_m$---excluding zero, which corresponds to the scalar representations. 
The decomposition for $C^*(T_\xi M)$ is obtained from the decomposition of its fibers, and it is important to observe how the representation theory of the Heisenberg group is reflected in the structure theory of the symbol algebra $C^*(T_\xi M)$.

Now that we understand how to decompose the symbol algebra $C^*(T_\xi M)$,
we can study its $K$-theory.
A quotient map is the algebraic analog of an inclusion of topological spaces.
And just as a pair of spaces $(X,Y)$ with $Y\subseteq X$ gives rise to a long exact sequence in cohomology, so does a short exact sequence of $C^*$-algebras.
In the present case, the part of the $K$-theory sequence that is relevant to our computation is,
\[ K^1(\xi^*) \to K^0(M\times \RR^{\times}) \to K_0(C^*(T_\xi M))\to 0.\] 
Let's analyze the three terms in this sequence is.
Starting from the left, since $\xi$ is a symplectic bundle (and hence $K$-orientable), the Thom isomorphism gives
\[K^1(\xi^*)\cong K^1(M).\]
From $\RR^{\times}\approx \RR\cup \RR$ we obtain
\[K^0(M\times \RR^{\times})\cong K^0(M\times \RR) \oplus K^0(M\times \RR) \cong K^1(M)\oplus K^1(M).\]
Finally, combining various isomorphisms,  
\[ K_0(C^*(T_\xi M)) \cong K^0(T^*M) \cong K^0(\xi\times \RR) \cong K^1(M).\]
Thus, our $K$-theory sequence can be expanded as follows,
\[ \xymatrix{    K^1(\xi^*) \ar[r]\ar[d]_{\cong} 
               & K^0(M\times \RR^{\times}) \ar[r]\ar[d]_{\cong} 
               & K_0(C^*(T_\xi M)) \ar[r]\ar[d]^{\cong} 
               & 0 \\
                 K^1(M) \ar[r] 
               & K^1(M)\oplus K^1(M) \ar[r] 
               & K^1(M) \ar[r]  
               & 0. }
\]
Notice that all the ingredients in this diagram are topological, except for the $C^*$-algebra $C^*(T_\xi M)$.
In order to compute the vertical arrow at the right of the diagram (the desired isomorphism $K_0((T_\xi M))\cong K^1(M)$) it suffices to work out explicitly what the other maps in the diagram are.
Then we will lift the element $\sigma_\xi(P) = \{P_m, m\in M\}$ from $K_0(C^*(T_\xi M)$ to $K^0(M\times \RR^{\times})$, and chase it through the diagram.
This is how we arrive at Theorem \ref{thm:contact} as a special case of Theorem \ref{thm:topind}.

Given the expository nature of this note we could stop here.
But for the reader who has read this far, it seems unfair not to give some final hints as to how the calculation is completed.
Glossing over many technical details (among other things, the complications arising from the fact that we are dealing with {\em unbounded} $K$-theory elements), we sketch the main idea.

Observe, before we proceed, that the vertical maps in the diagram (apart from the one on the right, which is the one we want to compute) are topologically well-understood.
What we need to identify are the two maps in the bottom row
by carefully analyzing the maps in the top row.
This translation from the top row into the bottom row is achieved by the techniques of {\em quantization}, or rather its reverse: ``passage to the classical limit''.

The result of this analysis is that we can identify the bottom-row map on the left as the diagonal embedding
\[ K^1(M) \to K^1(M)\oplus K^1(M)\;\colon\; a\mapsto a\oplus a.\]
This knowledge fixes the bottom-row map on the right,
\[ \xymatrix{    K^0(M\times \RR^{\times}) \ar[r]\ar[d]_{\cong} 
               & K_0(C^*(T_\xi M)) \ar[d]^{\cong} \\
                 K^1(M)\oplus K^1(M) \ar[r] 
               & K^1(M) }
\]
up to an automorphism of $K^1(M)$. It turns out to be the map
\[ K^1(M)\oplus K^1(M) \to K^1(M) \;\colon\; (a, b) \mapsto a\cdot b^{-1}.\]
Now we must lift the family $\{P_m, m\in M\}$ from $K_0(C^*(T_HM))$ to $K^0(M\times \RR^{\times})$. 
Here we must remember the explicit connection between the decomposition of $C^*(T_HM)$ and the representation theory of the Heisenberg group.
Very roughly (I am cheating a bit here),
the Heisenberg symbol lifts to the family of operators $P_{m,t}$ parametrized by $(m,t)\in M\times \RR^{\times}$, where $P_{m,t} = \pi_t(P_m)$ for $t>0$, while  $P_{m,-t}=\pi_{t}(P_m^{op})$.
Chasing this element through the diagram, we find
\[ K^0(M\times \RR^{\times})\to K^1(M)\oplus K^1(M) \;\colon\; \{P_{m,t}\} \mapsto (P_{m,+1}, P_{m, -1}) = (\pi_{+1}(P_m), \pi_{+1}(P_m^{op})),\]
and then
\[ K^1(M)\oplus K^1(M)\to K^1(M) \;\colon\; (\pi_{+1}(P_m), \pi_{+1}(P_m^{op})) \mapsto \pi_{+1}(P_m)\pi_{+1}(P_m^{op})^{-1}.\]
The result is an explicit formula for the isomorphism between the $K$-theories of the noncommutative symbol algebra $C^*(T_\xi M)$ and the topological space $T^*M$.
 
\begin{proposition}
Let $P$ be a $\xi$-elliptic operator on a closed contact manifold $(M, \xi)$.
With the canonical identifications
\[ K_0(C^*(T_\xi M)) \cong K^0(T^*M) \cong K^1(M),\]
we have
\[ [\{P_m, m\in M\}] = [\pi_{+1}(P_m)\pi_{+1}(P_m^{op})^{-1}].\]
\end{proposition}
This computation in $K$-theory is how we derive Theorem \ref{thm:contact}
from Theorem \ref{thm:topind}.
In turn, as we have seen, the World's Simplest Index Theorem \ref{thm:simple} is a special case of Theorem \ref{thm:contact}.
We hope that this discussion of the solution of this problem reveals something of the power and role of noncommutative topology as a bridge between classical analysis and topology.

\bibliographystyle{amsxport}
%

\begin{bibdiv}
\begin{biblist}

\bib{Ho71}{incollection}{
      author={H{\"o}rmander, Lars},
       title={On the index of pseudodifferential operators},
        date={1971},
   booktitle={Elliptische {D}ifferentialgleichungen, {B}and {II}},
   publisher={Akademie-Verlag},
     address={Berlin},
       pages={127\ndash 146. Schriftenreihe Inst. Math. Deutsch. Akad.
  Wissensch. Berlin, Reihe A, Heft 8},
}

\bib{Me97}{incollection}{
      author={Melrose, Richard},
       title={Homology and the {H}eisenberg algebra},
        date={1997},
   booktitle={S\'eminaire sur les \'{E}quations aux {D}\'eriv\'ees
  {P}artielles, 1996--1997},
   publisher={\'Ecole Polytech.},
     address={Palaiseau},
       pages={Exp.\ No.\ XII, 11},
        note={Joint work with C. Epstein and G. Mendoza},
}

\bib{Ep04}{incollection}{
      author={Epstein, Charles~L.},
       title={Lectures on indices and relative indices on contact and
  {CR}-manifolds},
        date={2004},
   booktitle={Woods {H}ole mathematics},
      series={Ser. Knots Everything},
      volume={34},
   publisher={World Sci. Publ., Hackensack, NJ},
       pages={27\ndash 93},
}

\bib{Bo79}{article}{
      author={Boutet~de Monvel, Louis},
       title={On the index of {T}oeplitz operators of several complex
  variables},
        date={1978/79},
        ISSN={0020-9910},
     journal={Invent. Math.},
      volume={50},
      number={3},
       pages={249\ndash 272},
         url={http://dx.doi.org/10.1007/BF01410080},
}

\bib{vEc}{article}{
      author={van Erp, Erik},
       title={The index of hypoelliptic operators on foliated manifolds},
        date={2011},
     journal={J. Noncommut. Geom.},
       pages={in press},
}

\bib{Ho67}{article}{
      author={H{\"o}rmander, Lars},
       title={Hypoelliptic second order differential equations},
        date={1967},
        ISSN={0001-5962},
     journal={Acta Math.},
      volume={119},
       pages={147\ndash 171},
}

\bib{FS74}{article}{
      author={Folland, G.~B.},
      author={Stein, E.~M.},
       title={Estimates for the {$\bar \partial _{b}$} complex and analysis on
  the {H}eisenberg group},
        date={1974},
        ISSN={0010-3640},
     journal={Comm. Pure Appl. Math.},
      volume={27},
       pages={429\ndash 522},
}

\bib{Ro78}{article}{
      author={Rockland, Charles},
       title={Hypoellipticity on the {H}eisenberg group:
  representation-theoretic criteria},
        date={1978},
        ISSN={0002-9947},
     journal={Trans. Amer. Math. Soc.},
      volume={240},
       pages={1\ndash 52},
}

\bib{BG88}{book}{
      author={Beals, Richard},
      author={Greiner, Peter},
       title={Calculus on {H}eisenberg manifolds},
      series={Annals of Mathematics Studies},
   publisher={Princeton University Press},
     address={Princeton, NJ},
        date={1988},
      volume={119},
}

\bib{Ta84}{article}{
      author={Taylor, Michael~E.},
       title={Noncommutative microlocal analysis. {I}},
        date={1984},
        ISSN={0065-9266},
     journal={Mem. Amer. Math. Soc.},
      volume={52},
      number={313},
       pages={iv+182},
}

\bib{vEb}{article}{
      author={van Erp, Erik},
       title={The {A}tiyah-{S}inger formula for subelliptic operators on a
  contact manifold, {P}art {II}},
        date={2010},
     journal={Ann. of Math.},
       pages={in press},
}

\end{biblist}
\end{bibdiv}

\end{document}